\input amstex.tex
\documentstyle{amsppt}
\overfullrule=0pt
\leftheadtext{C.D.Hill and M.Nacinovich}


\magnification=1200
\vsize=8.5truein

\font\sc=cmcsc10

\hyphenation{pseu-do-con-ca-ve}
\hyphenation{s-o-n-d-e-r-f-o-r-s-c-h-u-n-g-s-b-e-r-e-i-c-h}
\newcount\q
\newcount\x
\newcount\t

\long\def\se#1{\advance\q by 1
\x=0  \t=0 \bigskip
\noindent
\S\number\q \quad
{\bf {#1}}\par
\nopagebreak}

\long\def\thm#1{\advance\x by 1
\bigskip\noindent%
{\sc Theorem \number\q.\number\x}
\quad{\sl #1} \smallskip\noindent}

\long\def\prop#1{\advance\x by 1
\bigskip\noindent%
{\sc Proposition \number\q.\number\x}
\quad{\sl #1} \smallskip\noindent}

\long\def\thml#1#2{\advance\x by 1
\bigskip\noindent
{\sc Theorem \number\q.\number\x ({#1})}
\quad{\sl #2} \smallskip\noindent}

\long\def\lem#1{\advance\x by 1
\medskip\noindent
{\sc Lemma \number\q.\number\x}
\quad{\sl #1} \smallskip\noindent}

\def\supp{\roman{supp}}

\long\def\form#1{\global\advance\t by 1
$${#1} \tag \number\q.\number\t$$}
\long\def\cor#1{\advance\x by 1
\bigskip\noindent%
{\sc Corollary \number\q.\number\x}%
\quad{\sl #1} \smallskip\noindent}

\def\dimo{\smallskip\noindent {\sc Proof}\quad}

\topmatter
\title
Weak pseudoconcavity and the maximum modulus principle
\endtitle
\author
C.Denson Hill and Mauro Nacinovich
\endauthor
\address C.Denson Hill - Department of Mathematics, SUNY at Stony Brook,
Stony Brook NY 11794, USA \endaddress
\email dhill\@math.sunysb.edu \endemail
\address Mauro Nacinovich - Dipartimento di Matematica "L.Tonelli" -
via F.Buonarroti, 2 - 56127 PISA, Italy \endaddress
\email nacinovi\@ dm.unipi.it\endemail
\endtopmatter
\rightheadtext{Weak pseudoconcavity ...}
\document
In this paper we focus on the maximum modulus principle
and weak unique continuation for $CR$ functions
on an abstract {\it almost} $CR$ manifold $M$. It is known that some
assumption must be made on $M$ in order to have either of these:
it suffices to consider the standard $CR$ structure on the 
sphere $S^3$ in $\Bbb C^2$ to see that
the maximum modulus principle is not valid in the presence of strict
pseudoconvexity. For weak unique continuation, Rosay [R] has shown
by an example that there is a strictly pseudoconvex $CR$ structure on
$\Bbb R^3$, which is a perturbation of the aforementioned
standard $CR$ structure on $S^3$, such that there exists a smooth
$CR$ function $u$, $u\not\equiv 0$, with $u\equiv 0$ on a nonempty open
set. However positive results were obtained 
in [DCN] under the assumption of {\it pseudoconcavity}
and in [HN] under the assumption
of {\it essential pseudoconcavity} (and also finite kind for the maximum
modulus principle).\par
Here we investigate these matters under the assumption of {\it weak
pseudoconcavity} on $M$, which is a more general notion than that of
essential pseudoconcavity, insofar as it drops the minimality (and the
finite kind) hypothesis on $M$. We obtain sharp results involving
propagation along Sussmann leaves. The core of our argument is that
on a weakly pseudoconcave $M$ the square of the modulus of a $CR$
function is subharmonic with respect to a degenerate-elliptic
operator $P$ on $M$. We employ a maximum principle 
for real valued functions which is in the
spirit of [Hf], [Ni], [B], [H].\par
In order to understand our motivation in considering the weak
pseudoconcavity condition on $M$, the reader is referred to the
examples in [HN].
\se{Weak pseudoconcavity of almost $CR$ manifolds}
An abstract smooth almost $CR$ manifold
of type $(n,k)$ consists of: a connected smooth paracompact manifold
$M$ of dimension $2n+k$, a smooth subbundle $HM$ of $TM$ of rank
$2n$, and a smooth complex structure $J$ on the fibers of $HM$.\par Let
$T^{0,1}M$ be the complex subbundle of the complexification $\Bbb CHM$
of $HM$, which corresponds to the $-\sqrt{-1}$ eigenspace of $J$:
\form{T^{0,1}M\,=\,\{X+\sqrt{-1}JX\, \big{|}\, X\in HM\}\, .} We say
that $M$ is a $CR$ manifold if, moreover, the formal integrability
condition \form{\left[ \Cal C^\infty(M,T^{0,1}M),\,
\Cal C^\infty(M,T^{0,1}M)\right]\,
\subset\, \Cal C^\infty(M,T^{0,1}M)\, } 
\noindent
holds.

\par
Next we define ${T^*}^{1,0}M$ as the
annihilator of $T^{0,1}M$ in the complexified cotangent bundle $\Bbb C
T^*M$.  We denote by $Q^{0,1}M$ the quotient bundle $\Bbb C T^* M /
{T^*}^{1,0}M$, with projection $\pi_Q$.  It is a rank $n$ complex vector
bundle on $M$, dual to $T^{0,1}M$.  The $\bar\partial_M$--operator
acting on smooth functions is defined by $\bar\partial_M= \pi_Q\circ d$.
A local trivialization of the bundle $Q^{0,1}M$ on an open set $U$ in
$M$ defines $n$ smooth sections $\bar L_1$, $\bar L_2$, $\hdots$, $\bar
L_n$ of $T^{0,1}M$ in $U$; hence
\form{\bar\partial_M u\, = \,
\left(\bar L_1u, \bar L_2u,\hdots, \bar L_nu\right)\, ,}
where $u$ is a

function in $U$.  Solutions $u$ of $\bar\partial_Mu=0$ are called $CR$
functions.
\par
The {\it characteristic bundle} $H^0M$ is defined to be the
annihilator of $HM$ in $T^*M$. Its purpose it to parametrize
the Levi form: recall that the {\it Levi form}
of $M$ at $x$ is defined for $\xi\in H^0_xM$ and $X\in H_xM$ by
\form{\Cal L(\xi;X)\, = \, d\tilde\xi(X,JX)=\langle\xi,
[J\tilde X,\tilde X]\rangle\, ,}
where $\tilde\xi\in\Cal C^\infty(M,H^0M)$ and $\tilde
X\in\Cal C^\infty(M,HM)$
are smooth extensions of $\xi$ and $X$.
For each fixed $\xi$ it is a Hermitian quadratic form for the
complex structure $J_x$ on $H_xM$.
\smallskip
Denote by $H^{1,1}M$ the smooth subbundle of the tensor bundle
$HM\otimes_M HM$ whose fiber $H^{1,1}_xM$ at $x\in M$ is the real
vector subspace of $H_xM\otimes H_xM$ generated by the tensors
of the form $v\otimes v+(Jv)\otimes(Jv)$ for $v\in H_xM$.
$H^{1,1}M$ is the bundle of {\sl Hermitian symmetric tensors}
in $HM\otimes_MHM$. For each $x\in M$ and $\xi\in H^0M$ the
Levi form $\Cal L(\xi,\,\cdot\,)$ defines a linear form
$\Cal L_\xi:H^{1,1}M@>>>\Bbb R$ such that
\form{\Cal L_{\xi}\left(v\otimes v+ (Jv)\otimes (Jv)\right)
=\Cal L(\xi,v)
\qquad \forall v\in H_xM\, .}
For $x\in M$ let us denote by
$\bar\Gamma H^{1,1}_xM$ the convex hull of
$\{v\otimes v+(Jv)\otimes (Jv)\, |\, v\in H_xM\}$ and by
$\Gamma H^{1,1}M$ its interior (in $H^{1,1}_xM\simeq\Bbb R^{n^2}$).
They are the closed cone of nonnegative Hermitian symmetric tensors and
the open cone of positive Hermitian symmetric tensors
of $H_xM\otimes H_xM$, respectively.
The disjoint union
$\Gamma H^{1,1}M=\cup_{x\in M}{\Gamma H^{1,1}_xM}$
is an open subset of $H^{1,1}M$ and the restriction of the projection
onto the base:
\form{\pi:\Gamma H^{1,1}M @>>> M}
is a smooth fiber bundle, whose fibers are open convex cones in
$\Bbb R^{n^2}$. Note that the choice of a smooth Hermitian metric $h$ on
the fibers of $HM$ defines an exponential map
\form{\roman{exp}_h:H^{1,1}M @>>> \Gamma H^{1,1}M\, ,}
giving a smooth bundle isomorphism between $\Gamma H^{1,1}M$ and
$H^{1,1}M$.
\smallskip
\noindent
{\sc Definition}\quad We say that an abstract almost $CR$ manifold $M$
is {\it weakly pseudoconcave} iff
for every $x\in M$ there is an open neighborhood
$U$ of $x$ in $M$ and a smooth section
$\Omega\in\Cal C^\infty(U,\Gamma H^{1,1}M)$ such that
\form{\Cal L_\xi(\Omega)=0\qquad\forall x\in U,\;
\xi\in H^0_xM\, .}
\medskip
\noindent
{\sc Remark}\quad
Every abstract almost $CR$ manifold, whose Levi form vanishes
identically, is trivially weakly pseudoconcave. However, when
$k>0$, such a manifold is not necessarily {\it essentially pseudoconcave}
in the sense of Definition A of [HN].\par
An abstract almost $CR$ manifold of type $(n,0)$ is the same thing 
as an {\it almost complex manifold}; such manifold can be regarded
as being essentially pseudoconcave, and hence weakly pseudoconcave.
In this case the $CR$ functions will be called {\it almost holomorphic
functions}.
\smallskip
We shall need the following results from [HN]:
\prop{Let $M$ be an abstract almost $CR$ manifold of type $(n,k)$.
Then $M$ is weakly pseudoconcave if and only if there exists
a smooth Hermitian metric $\bold h$ on the fibers of $HM$ such that
\form{\roman{trace}_{\bold h}\left(\Cal L(\xi,\,\cdot\,)\right)\, = \, 0\,
,\qquad
\forall \xi\in H^0M\, .}}
\prop{Let $M$ be an abstract almost $CR$ manifold of type $(n,k)$.
If $M$ is weakly pseudoconcave then
\form{\cases
\text{For each $\xi\in H^0M$ the Levi form $\Cal L(\xi,\,\cdot\,)$ is
either $0$}\\
\text{or has at least one positive and one negative eigenvalue.}
\endcases}
If $\Cal D:=\Cal C^\infty(M,HM)+[C^\infty(M,HM),C^\infty(M,HM)]$ is a
distribution of constant rank, then (\number\q.\number\t) is also
sufficient for $M$ to be weakly pseudoconcave.}
\prop{Under the assumptions of Proposition 1.1, let $U$ be an open
subset of $M$ on which $X_1,\hdots,X_n\in\Cal C^\infty(U,HM)$
give at each point $y\in U$ an $\bold h$-orthonormal basis of
the complex Hermitian vector space $H_yM$. Set
$\bar L_j=X_j+iJX_j$
and $L_j=X_j-iJX_j$, for $j=1,\hdots, n$.
Then there are smooth complex valued functions $\beta^r$ ($1\leq r\leq n$)
on $U$ such that
\form{i\dsize\sum_{j=1}^n{[L_j,\bar L_j]}=
\dsize\sum_{r=1}^n{\left(\beta^r\,L_r+\bar\beta^r\,\bar L_r\right)}
\quad\text{in}\quad U\, .}}
\edef\formula1-3{\number\q.\number\t}
\edef\propo1-3{\number\q.\number\x}
Let $L=X-iJX$  be one of the $L_j$'s from Proposition \number\q.\number\x.
We have
\form{\matrix\format\r&\l\\
\Re L\bar L&\, = \, X^2+\left(JX\right)^2\\
\Im L\bar L&\, = \, [X,JX]\, .
\endmatrix}
Let $u$ be a $CR$ function in $U$, and consider $|u|^2=u\bar u$. Since
\form{\bar L\, |u|^2\, = \, \left(\bar L u\right)\bar u
+u\,\bar L\, \bar u\, ,}
and $\bar L\, u=0$, we obtain
\form{L\bar L\, |u|^2\, = \, \left| Lu\right|^2\, + \, u\,
\left[L,\bar L\right]\, \bar u
\, .}
It follows that
\form{\matrix\format\r&\l\\
\left(\dsize\sum_{j=1}^n{L_j\bar L_j}\right)\, |u|^2 &\, =\,
\dsize\sum_{j=1}^n{\left| L_ju\right|^2}\, +\,
u\left(\dsize\sum_{j=1}^n{\left[L_j,\bar L_j\right]}\right)\, \bar u\\
&\, =\,\dsize\sum{\left| L_ju\right|^2}\, + \, u\, \left(
\frac{1}{i}\dsize\sum_{r=1}^n{\bar\beta^r\bar L_r}\right)\,\bar u\\
&\, = \, \dsize\sum_{j=1}^n{\left| L_ju\right|^2}\, + \,
\frac{1}{i}\left(\dsize\sum_{r=1}^n{\bar\beta^r\bar L_r}\right)\, |u|^2\, ,
\endmatrix}
because of (\formula1-3). Hence
\form{\left\{
\Re\left(\dsize\sum_{j=1}^n{L_j\bar L_j}\right)+\Im \left(
\dsize\sum_{j=1}^n{\beta^jL_j}\right)\right\}\, |u|^2\, =\,
\dsize\sum_{j=1}^n{\left| L_ju\right|^2}\, \geq\, 0\, .}
A similar calculation shows that
\form{
\left\{
\Re
\left(
\dsize\sum_{j=1}^n{L_j\bar L_j}
\right)+
\Im
\left(
\dsize\sum_{j=1}^n{\beta^jL_j}
\right)
\right\}
\Re u\, = \, 0\, .}
\edef\lalli{\number\q.\number\t}
Let $P_U$ denote the {\sl real} operator inside the curly brackets.
It has the form
\form{\dsize\sum_{j=1}^n{\left(X_j^2+\left(JX_j\right)^2\right)}\, + \,
X_0\, ,}
\noindent
where the $X_1,\hdots, X_n,JX_1,\hdots, JX_n$ provide a basis
for $HM$ at each point of $U$, and $X_0\in\Cal C^\infty(U,HM)$.
\edef\formula1-10{\number\q.\number\t}
\prop{Let $M$ be a weakly pseudoconcave almost $CR$ manifold
of type $(n,k)$. Then one can construct a smooth real linear
second order partial differential operator $P$ on $M$ such that:
\roster
\item"($i$)" each $x_0\in m$ has a neighborhood $U$ in which $P$
can be written in the form (\formula1-10);
\item"($ii$)" if $u$ is a $\Cal C^2$ $CR$ function on $M$, then
$Pu=0$ and $P|u|^2\geq 0$ on $M$.
\endroster}
\dimo It suffices to take
\form{P\, = \, \dsize\sum_{U}{\psi_U\, P_U}\, ,}
where $\{\psi_U\}$ is a nonnegative partition of unity subordinate to
a covering $\{U\}$ of $M$ by open sets $U$, as in Proposition \propo1-3.
Indeed ($ii$) is then obvious, while ($i$) follows because
$P_U$ and $P_V$ have the same principal symbol on $U\cap V$.
\medskip
\edef\Propo1-4{\number\q.\number\x}
\se{Sussmann leaves}
In this section we collect the results which we shall need
concerning the Sussmann leaves of an arbitrary set $\Cal D$ of
smooth real vector fields on a smooth paracompact manifold $M$
of real dimension $N$. In our final application, $M$ will be an
abstract almost $CR$ manifold, and $\Cal D=\Cal C^\infty(M,HM)$.
However, in our discussion of the maximum principle for real valued
functions, in the next section, we shall be in this more general
situation.
\par
Let $x_0\in M$ and $\Omega$ be an open subset of $M$ containing $x_0$.
The {\it Sussmann leaf} $\Cal F(x_0,\Omega)$ of $\Cal D$ in $\Omega$
through $x_0$ is defined to be the set of points $x\in\Omega$
for which there exist finitely many smooth curves
$s_j:[0,1]@>>>\Omega$, for $j=1,\hdots,\ell$, such that:
\form{
\cases
\dot{s}_j(t)\in\Cal D_{s_j(t)}\qquad\text{for}\quad 0\leq r\leq 1\quad
\text{and}\quad j=1,2,\hdots,\ell\, ;\\
s_j(0)=x_0,\quad s_j(0)=s_{j-1}(1) \quad\text{for}\quad j=2,\hdots,\ell
\qquad\text{and}\quad s_\ell(1)=x\, .
\endcases}
Note that $\Cal F(x,\Omega)=\Cal F(x_0,\Omega)$ for all
$x\in\Cal F(x_0,\Omega)$. Sussmann proved in [S] that
$\Cal F(x_0,\Omega)$  is always a smooth immersed (but not necessarily
embedded) submanifold of $\Omega$. Note also that
$T_x\Cal F(x_0,\Omega)\supset\Cal D_x$ for all $x\in\Cal F(x_0,\Omega)$.
We say that $M$ is {\it minimal} at $x_0$ in $M$ iff for every open
neighborhood $U$ of $x_0$ in $M$, the Sussmann leaf
$\Cal F(x_0,U)$ contains an open neighborhood of $x_0$ in $M$.
The manifold $M$ is said to be {\it minimal} if it is minimal at
each point. This condition is equivalent to the nonexistence of a
lower dimensional smooth submanifold $S$ of $M$ with
$x_0\in S$ and $T_xS\supset\Cal D_x$ for every $x\in S$.
\smallskip
Next we recall the definition of the set $N_eF$ of
{\it exterior conormals} to a closed subset $F$ of $M$: it is the subset
of $T^*M$ consisting of all the nonzero $\xi_0\in T_{x_0}^*M$,
with $x_0\in F$, for which there exists a smooth real valued function
$f$ on $M$ with $df(x_0)=\xi_0$ and $f(x)\leq f(x_0)$ for all
$x\in F$.
\par
In what follows we shall use the well known trapping lemma 
(see for instance [Ho I, Theorem 8.5.11, p.304]):
\prop{Let $F$ be a closed subset of $M$. If
\form{\xi(X)=0\quad\text{for all}\quad \xi\in N_eF\quad\text{and}\quad
\text{all}\quad X\in\Cal D\, ,}
then $\Cal F(x,M)\subset F$ for every $x\in F$.}
\edef\propo2-1{\number\q.\number\x}
\edef\sezion2{\number\q}
\se{A maximum principle for real valued functions}
Let $M$ and $\Cal D$ be as in section \sezion2. We shall consider a
smooth real second order linear partial differential operator $P$ on
 $M$ with the following property: Given $x_0\in M$, there is an open
neighborhood $U$ of $x_0$ in $M$, and $Y_0,Y_1,\hdots,Y_\ell\in\Cal D$
such that
\form{\cases
Y_1,\hdots,Y_\ell\quad\text{generate $\Cal D$ in $U$}\, ,\\
P\,=\,\dsize\sum_{j=1}^\ell{Y_j^2}\, + \, Y_0\quad\text{in $U$}.
\endcases}
\edef\formu3-1{\number\q.\number\t}
\thm{Let $\Omega$ be an open subset of $M$, $x_0\in\Omega$,
$u\in\Cal C^2(\Omega,\Bbb R)$ and $Pu\geq 0$ along
 $\Cal F\left(x_0,\Omega\right)$. 
If $u(x)\leq u(x_0)$ for all $x\in\Cal F(x_0,\Omega)$,
then $u$ is constant along $\overline{\Cal F(x_0,\Omega)}\cap\Omega$.}
\dimo
For the proof we can, without loss of generality, assume that
$\Omega=M=\Cal F(x_0,\Omega)$ and $Pu\geq 0$ on $M$. \par
Let $F$ denote the closed subset
$\{x\in M\, | \, u(x)=u(x_0)\,\}$. We want to show that $F=M$.
Assume by contradiction that $F\neq M$; i.e., that $F$ does not
contain $\Cal F(x_0,M)$. By Proposition \propo2-1\; there exist
$x_1\in\partial F$, $\xi\in T_{x_1}^*M$ with $\xi\in N_eF$ and
$Y\in\Cal D$ such that $\xi(Y)\neq 0$. This implies the following:
there is a coordinate patch $U\simeq\{y\in\Bbb R^N\, | \, |y|<R\}$
containing $x_1$, with $0<|y(x_1)|=r<R$, such that
\roster
\item"($i$)" $P=\dsize\sum_{j=1}^\ell{Y_j^2}+Y_0$ in $U$ with
$Y_0,Y_1,\hdots,Y_\ell\in\Cal D$;
\item"($ii$)" $Y_{j_0}(|y|^2)\neq 0$ at $x_1$ for some $j_0$ with
$1\leq j_0\leq \ell$;
\item"($iii$)" $u(x)<u(x_0)=u(x_1)$ if $x\in U$ and
$|y(x)|\leq r$, $x\neq x_1$.
\endroster
Let $\gamma>0$. Then
\form{P\left(e^{-\gamma |y|^2}\right)
\, = \,e^{-\gamma |y|^2}
\left\{
\gamma^2 \dsize\sum_{j=1}^\ell{\left|Y_j\left(|y|^2\right)\right|^2}\, + \,
O(\gamma)\,\right\}}
is positive on a neighborhood of $x_1$ for $\gamma>0$ sufficiently large.
Fix $\gamma>0$ and $\epsilon>0$ in such a way that
$0<\epsilon<R-r$ and $P(\exp(-\gamma |y|^2))>0$ when $x\in U$ and
$|y(x)-y(x_1)|\leq\epsilon$. For $\delta>0$ set
$v_\delta=u+\delta\left(e^{-\gamma|y|^2}-e^{-\gamma r^2}\right)$.
Then $Pv_\delta>0$ for $|y(x)-y(x_1)|\leq\epsilon$. Note that
$v_\delta(x)<u(x)$ when $|y(x)|>r$. On the other hand, $u(x)<u(x_0)$
if $|y(x)|\leq r$ and $|y(x)-y(x_1)|=\epsilon$. Thus for $\delta>0$
sufficiently small, we obtain that $v_\delta(x)<u(x_0)=u(x_1)$ on the
boundary of $\omega=\{x\in U\, | \, |y(x)-y(x_1)|<\epsilon\}$.
Since $v_\delta(x_1)=u(x_1)=u(x_0)$, the restriction of $v_\delta$
to $\overline{\omega}$ has a maximum at some point $x_2\in\omega$.
But at $x_2$ we would then have that $Pv_\delta(x_2)\leq 0$,
which contradicts the inequality $Pv_\delta>0$ we have established
in $\omega$. Thus $F=M$ and the theorem is proved, after using
continuity of $u$ to pass to the closure of the Sussmann leaf.
\edef\trema{\number\q.\number\x}
\se{Weak unique continuation}
In this section we return to a smooth manifold $M$ which is an abstract
{\it almost} $CR$ manifold of type $(n,k)$, and $\Cal D$ will be
$\Cal C^\infty(M,HM)$. In this situation, for any open $\Omega\subset M$
and $x_0\in M$, the Sussmann leaf $\Cal F(x_o,\Omega)$ is 
itself a smooth
abstract almost $CR$ manifold of type $(n,h)$ for some $h\leq k$. \par
The next theorem is an improvement of the weak unique continuation
result of [DCN, Theorem 4.1], [HN, Theorem 5.1].
\thm{Assume that $M$ is weakly pseudoconcave. Let
$u\in L^2_{\text{loc}}(M)$ satisfy the following:
\form{\boxed{\matrix\format\l\\
\text{for every $\bar L\in\Cal C^\infty(M,T^{0,1}M)$,
$\bar Lu\in L_{\text{loc}}^2(M)$}\\
\text{and there exists $\kappa_{\bar L}\in L^\infty_{\text{loc}}(M)$
such that}\\
\quad \left|\bar L u(x)\right|\leq \kappa_{\bar L}(x)\, |u(x)|\quad
\text{a.e. in}\quad M\, .
\endmatrix}}
Then $\Cal F(x,M)\subset\supp\, u$ for every $x\in\supp\, u$.}
\edef\formula4-1{\number\q.\number\t}
\dimo
We use again Proposition \propo2-1. Indeed under the contrary assumption,
there exists a $\xi\in N_e(\supp\, u)$ such that $\xi(X)\neq 0$ for
some $X\in HM$. We obtain a contradiction by using the Carleman type
estimate given by the following theorem.
\thm{Let $M$ be a weakly pseudoconcave abstract almost $CR$ manifold
of type $(n,k)$. Let $\phi$ be a real valued smooth function on $M$
and $x_0\in M$ a point where $\phi(x_0)=0$ and
$d\phi(x_0)\notin H^0M$. Then we can find $A>0$, $C>0$, $\tau_0>0$
and an open neighborhood $U$ of $x_0$ in $M$ such that:
\form{\matrix\format\l\\
\sqrt{\tau}\cdot\| f\cdot\exp(\tau(\phi+A\phi^2))\|_0
\leq c \|\bar\partial_M f\cdot\exp(\tau(\phi+A\phi^2))\|_0\\
\qquad\qquad
\forall f\in\Cal C^\infty_0(U)\, ,\quad\forall\tau\geq\tau_0\, .
\endmatrix}}
Here the $L^2$-norms $\|\,\cdot\,\|_0$ are computed using any smooth
Riemannian metric on $M$ and any smooth Hermitian metric on the
fibers of $Q^{0,1}M$.
\par
\edef\teorema4-1{\number\q.\number\x}
Theorem \number\q.\number\x \; is just Theorem 5.2 of [HN], with
"weakly pseudoconcave" replacing "essentially pseudoconcave" in the
hypothesis. In fact the proof of Theorem 5.2 in [HN] does not use
the minimality assumption on $M$, which is part of the definition
of essential pseudoconcavity, but only uses the weak pseudoconcavity.
\cor{Assume that $M$ is weakly pseudoconcave. Let $u$
be a continuous $CR$ function on $M$, and $x_0\in M$. Let $\omega$
be an open neighborhood of $x_0$ in $M$. If $u\equiv 0$ on
$\Cal F(x_0,M)\cap \omega$, then $u\equiv 0$ along $\Cal F(x_0,M)$.}
\dimo
We obtain the Corollary from Theorem \teorema4-1, after replacing
$M$ by $\Cal F(x_0,M)$.
\edef\corvo{\number\q.\number\x}
\cor{Let $M$ be a weakly pseudoconcave smooth abstract $CR$ manifold
of type $(n,k)$. Let $\frak L @>p>> M$ be a smooth complex $CR$ line
bundle over $M$, and $u$ be a continuous $CR$ section of $\frak L$
over $M$. If $x_0\in M$ and $\omega$ is an open neighborhood of
$x_0$ such that $u\equiv 0$ on $\Cal F(x_0,M)\cap \omega$, then
$u\equiv 0$ along $\Cal F(x_0,M)$.}
\dimo For the notion of a complex $CR$ line bundle we refer to
section 7 of [HN]. The corollary follows from Theorem \teorema4-1\;
because, according to formula (7.4) in [HN], the representative of
the section $u$, in any smooth (not necessarily $CR$) local
trivialization of $\frak L$, satisfies (\formula4-1).
\se{The maximum modulus principle}
In this section we have: $M$ is a smooth abstract almost $CR$ manifold
of type $(n,k)$, $\Omega$ is an open subset of $M$, and
$\Cal D=\Cal C^\infty(M,HM)$. Fix a point $x_0\in\Omega$ and set
$\Cal F=\Cal F(x_0,\Omega)$.
\lem{Let $u\in\Cal C^1(\Omega)$ be a $CR$ function in $\Omega$. Assume
that $u\left|_{\dsize\Cal F}\right.$ has values which lie along a
piecewise $\Cal C^1$-regular curve in $\Bbb C$. Then $u(x)=u(x_0)$ for
every $x\in\overline{\Cal F}$.}
\dimo It suffices to show that $u$ is locally constant along $\Cal F$,
and we can also assume that the values of $u$ lie on a $\Cal C^1$-regular
curve in $\Bbb C$. Let $\gamma$ be the $\Cal C^1$-regular curve in
$\Bbb C_z=\Bbb R_x\times\Bbb R_y$. Let $p_0\in\Cal F$ and $\omega$
be a connected open neighborhood of $u(p_0)$ in $\gamma$. If we take
$\omega$ sufficiently small, then there is an open neighborhood $\Omega$
of $u(p_0)$ in $\Bbb C$, and a real valued $\Cal C^1$ function $F(x,y)$
in $\Omega$ such that
\form{\omega\, = \, \{x+iy\in\Omega\, | \, F(x,y)=0\}\, ,\qquad
dF\neq 0 \quad\text{in}\quad \Omega\, .}
Choose a connected open neighborhood $V$ of $p_0$ in $\Cal F$
such that $u(V)\subset\omega$. Then $F(\Re u,\Im u)=0$ on $V$, so
\form{\matrix\format\r&\l\\
0&=\bar\partial_{\Cal F}F\, = \, F_u\bar\partial_{\Cal F}u+
F_{\bar u}\bar\partial_{\Cal F}\bar u\\
&=F_{\bar u}\bar\partial_{\Cal F}\bar u\qquad\text{and}\qquad
F_{\bar u}\neq 0\, ;\endmatrix}
hence $Xu=0$ in $V$ for every $X\in\Cal D$. This in turn implies that
$u$ is constant along $\Cal F$ in $V$, and hence along $\Cal F$.
\edef\lemma5-1{\number\q.\number\x}
\medskip
\noindent
{\sc Remark}\quad The lemma remains valid if we assume
 $u\in\Cal C^1(\Cal F)$ and
 $u$ is $CR$ on the almost $CR$ manifold $\Cal F$.
\thm{Let $M$ be a smooth abstract weakly 
pseudoconcave almost $CR$ manifold of
type $(n,k)$. Consider an open subset $\Omega$ of $M$ and a point
$x_0\in \Omega$. Let $u\in\Cal C^2\left(\Cal F(x_0,\Omega)\right)$
be a $CR$ function on the almost 
 $CR$ manifold $\Cal F(x_0,\Omega)$. Assume that
\form{\left| u(x_0)\right| \, = \, \sup_{\Cal F(x_0,\Omega)}{\left| u\right|}
\, .}
Then $u$ is constant along $\overline{\Cal F(x_0,\Omega)}\cap\Omega$.}
\dimo We observe that $\Cal F(x_0,\Omega)$ is a smooth abstract
almost
$CR$ manifold of type $(n,k)$ for some $h\leq k$. By Proposition \Propo1-4
there is a smooth real linear second order operator $P$ on
$\Cal F(x_0,\Omega)$ of the form (\formu3-1) such that $P|u|^2\geq 0$.
By Theorem \trema\; the real valued function $|u|^2$ is constant
along $\Cal F(x_0,\Omega)$. According to Lemma \lemma5-1, $u$ is constant
along $\Cal F(x_0,\Omega)$.
\edef\tmma{\number\q.\number\x}
\thm{Let $M$ be a smooth abstract weakly 
pseudoconcave almost $CR$ manifold of
type $(n,k)$. Consider a nonempty open subset $\Omega$ of $M$ and a point
$x_0\in\Omega$. Let $u\in\Cal C^2\left(\Cal F(x_0,\Omega)\right)$
be a $CR$ function on the almost
 $CR$ manifold $\Cal F(x_0,\Omega)$. Assume that
$M$ is minimal at $x_0$ and that $|u|$ has a local
weak maximum at $x_0$. Then $u$ is constant along
$\overline{\Cal F(x_0,\Omega)}\cap\Omega$.}
\dimo By our assumption $\Cal F(x_0,\Omega)$ is an open neighborhood of
$x_0$ in $\Omega$. Hence there is an open subset $\omega$ of $\Omega$,
containing $x_0$, such that
\form{\left| u(x_0)\right|\, = \, \sup_{\omega}{|u|}\, .}
By Theorem \tmma\; it follows that $u$ is constant along
$\Cal F(x_0,\omega)$, which is a neighborhood of $x_0$ in $\Omega$.
Corollary \corvo\; then implies that the function $u-u(x_0)$ is identically
zero along $\Cal F(x_0,\Omega)$.
\edef\tmmb{\number\q.\number\x}
\medskip
Recall that the notion of {\it essential pseudoconcavity} in [HN]
is weak pseudoconcavity plus minimality. Thus we obtain the following
improvement of Theorem 6.4 in [HM]:
\cor{Assume that $M$ is a smooth connected essentially pseudoconcave
abstract almost $CR$ manifold of type $(n,k)$. Let $u\in\Cal C^2(M)$
be a $CR$ function on $M$. If $|u|$ has a weak local maximum
at some point $x_0$ of $M$, then $u$ is constant on $M$.}
\edef\tmmc{\number\q.\number\x}
\medskip
\noindent
{\sc Remark 1}\quad In the statement of Theorem \tmma, Theorem \tmmb,
and Corollary \tmmc\; one can substitute $\Re u$ in place of
$|u|$, because of (\lalli). In particular if $M$ is as in
Corollary \tmmc, a $\Cal C^2$ $CR$ function on $M$, which is real valued
on a neighborhood of a point of $M$, is constant on $M$.
\smallskip
\noindent
{\sc Remark 2}\quad Suppose $M$ is an almost complex manifold. Then, 
according to Corollaries 4.3, 4.4, 5.4, the almost holomorphic functions
on $M$ obey weak unique continuation, and enjoy the usual form of
the maximum modulus principle. However in this situation the almost
holomorphic functions obey {\it strong} unique continuation, because of
(1.17), according to Theorem 17.2.6 in [Ho III].

\enddocument